\documentclass[11pt]{amsart}

\usepackage{amsmath,amssymb,amsfonts}
\usepackage[a4paper,margin=1.15in]{geometry}

\newenvironment{romanlist}{\begin{enumerate}}{\end{enumerate}}
\newenvironment{alphalist}{\begin{enumerate}}{\end{enumerate}}
\newenvironment{arabiclist}{\begin{enumerate}}{\end{enumerate}}

\theoremstyle{plain}
\newtheorem{theorem}{Theorem}[section]
\newtheorem{proposition}[theorem]{Proposition}
\newtheorem{lemma}[theorem]{Lemma}
\newtheorem{corollary}[theorem]{Corollary}
\theoremstyle{definition}
\newtheorem{remark}[theorem]{Remark}
\newtheorem{example}[theorem]{Example}

\numberwithin{equation}{section}
\newcommand{\kk}{K}
\newcommand{\Hess}{\operatorname{Hess}}
\newcommand{\adj}{\operatorname{adj}}
\newcommand{\im}{\operatorname{im}}
\newcommand{\ip}[1]{\langle #1\rangle}
\newcommand{\grad}{\nabla}
\newcommand{\JC}{\mathrm{JC}}
\newcommand{\HC}{\mathrm{HC}}

\begin{document}

\title[The Legendre potential and gradient Keller maps]{The Legendre potential and invertibility criteria for gradient Keller maps}

\author{El\.zbieta Adamus}
\address{Faculty of Applied Mathematics, AGH University of Krakow, al.~Mickiewicza 30, 30-059 Krak\'ow, Poland}
\email{esowa@agh.edu.pl}

\author{Zbigniew Hajto}
\address{Faculty of Mathematics and Computer Science, Jagiellonian University, ul.~\L{}ojasiewicza 6, 30-348 Krak\'ow, Poland}
\email{zbigniew.hajto@uj.edu.pl}

\subjclass[2020]{Primary 14R15; Secondary 12H05, 13N15, 37J30}

\keywords{Jacobian conjecture, Hessian conjecture, Nambu derivations, Legendre transform, Picard--Vessiot theory, Pascal finite automorphisms, cotangent lift}

\begin{abstract}
Let $h$ be a polynomial in $n$ variables over an algebraically closed field of characteristic zero whose Hessian determinant is a nonzero constant, let $h^{*}$ be its formal Legendre transform and $g=\ip{X,\grad h}-h$. The Nambu derivations of $\grad h$ satisfy $\delta_i(g)=X_i$ and $\delta^\alpha(g)=(\partial^\alpha h^{*})\circ\grad h$. Consequently, $\grad h$ is a polynomial automorphism if and only if $g$ belongs to the algebra generated by $\partial_1h,\dots,\partial_nh$, and then $g=h^{*}(\partial_1h,\dots,\partial_nh)$. Moreover, $\grad h$ is a polynomial automorphism if and only if $\delta^\alpha(g)=0$ for all multi-indices $\alpha$ with $|\alpha|=(\deg h-1)^{n-1}+2$. In three variables the order $2\deg h-1$ suffices and is optimal. In the gradient case the differential Galois criterion of Adamus, Crespo and Hajto thus reduces to a single condition. If the Hessian determinant is a nonzero polynomial, not necessarily constant, then $g$ is a primitive element of the field extension generated by the gradient. After an affine normalization of the gradient, $h^{*}$ is the sum of a Pascal difference series of $g$, and we give a gradient automorphism in four variables whose normal form is not Pascal finite. The results do not settle the open cases of the Jacobian conjecture in two variables and of the Hessian conjecture in four variables.
\end{abstract}

\maketitle

\section{Introduction}\label{sec:intro}

Throughout, $\kk$ is an algebraically closed field of characteristic zero and $\kk[X]=\kk[X_1,\dots,X_n]$. A polynomial map $F\colon\kk^n\to\kk^n$ is a Keller map if $\det JF\in\kk^{\times}$. The Jacobian conjecture $\JC_n$ asserts that every Keller map $\kk^n\to\kk^n$ has a polynomial inverse \cite{Kel39,BCW82,vdE00}. A polynomial $h\in\kk[X]$ is called a Keller potential if its Hessian determinant is a nonzero constant. The Hessian conjecture $\HC_n$ of Meng \cite{Men06} asserts that the formal Legendre transform of every Keller potential in $n$ variables is a polynomial. The implications $\JC_n\Rightarrow\HC_n$ and $\HC_{2n}\Rightarrow\JC_n$ are elementary \cite{Men06,MY26}.\footnote{%
In July 2026 L.~Alp\"oge announced an explicit Keller map $\kk^3\to\kk^3$ which is not injective, see \cite[Section~1]{MY26}. Hence $\JC_n$ fails for $n\ge3$. Meng and Yang \cite{MY26} deduced that $\HC_n$ fails for $n\ge5$, that $\HC_2$ holds by \cite{Dil91} and that $\HC_3$ holds by \cite{dB15}. The cases $\JC_2$ and $\HC_4$ remain open, and $\HC_4$ implies $\JC_2$. For four-variable potentials of degree at most four $\HC_4$ is proved in \cite{Ni26}, and the case of Lorentzian signature and degree at most five is treated in \cite{Liu26}.}

Two earlier lines of work meet in this paper. The first concerns the Nambu derivations of a Keller map and the differential Galois criterion for invertibility. Crespo and the second author \cite{CH11} introduced the Nambu derivations $\delta_1,\dots,\delta_n$ of a Keller map $F$, characterized by $\delta_i(F_j)=\delta_{ij}$, and proved that $F$ is invertible if and only if $\kk(X)$ with these derivations is a Picard--Vessiot extension of $\kk$. In \cite{ABH17} this was turned into an effective criterion through generalized Wronskians, tested by the algebra membership algorithm of \cite{SS88}, and in \cite{ACH19} into the surjectivity of a map $\varphi\colon \kk(X)\otimes_{\kk}E\to \kk(X)\otimes_{\kk}R$, which reduces to the $n$ conditions $1\otimes X_i\in\im\varphi$. The second line \cite{ABCH17,ABCH19,AB22,AH27} studies formal inversion by means of the Pascal difference operator $\Delta_F(P)=P\circ F-P$. Automorphisms $F$ with $\Delta_F^m(X)=0$ for some $m$ are called Pascal finite \cite{ABCH19}. They include the maps $X+H$ with strongly nilpotent $JH$ \cite[Theorem~5.1]{AH27}, and the alternating Pascal series computes the formal inverse of any $F=X+H$ \cite{Ad26}. In the gradient case $F=\grad h$ both approaches apply to one polynomial, the Legendre potential $g=\ip{X,\grad h}-h$. The main results are the following.

\begin{theorem}[Legendre criterion]\label{thm:intro-member}
Let $h\in\kk[X]$ be such that $\det\Hess h$ is a nonzero polynomial, $F=\grad h$ and $g=\ip{X,\grad h}-h$. Then the Nambu derivations of $F$ satisfy $\delta_i(g)=X_i$ for $1\le i\le n$, and $\kk(X)=\kk(F)(g)$. If $h$ is a Keller potential, the following are equivalent.
\begin{alphalist}
\item $F$ is a polynomial automorphism.
\item $g\in\kk[F_1,\dots,F_n]$.
\item $g\in\kk(F_1,\dots,F_n)$.
\end{alphalist}
If these hold, then $g=h^{*}(F_1,\dots,F_n)$, where $h^{*}$ is the formal Legendre transform of $h$.
\end{theorem}

\begin{theorem}[Nambu tower of the Legendre potential]\label{thm:intro-tower}
Let $h$ be a Keller potential of degree $d\ge2$, $N=(d-1)^{n-1}+2$, $x_0\in\kk^n$, and let $h^{*}$ be the formal Legendre transform of $h$ at $x_0$. For every multi-index $\alpha$ one has $\delta^\alpha(g)=(\partial^\alpha h^{*})\circ\grad h$ in $\kk[[X-x_0]]$. Consequently, $\grad h$ is a polynomial automorphism if and only if $\delta^\alpha(g)=0$ for all multi-indices $\alpha$ with $|\alpha|=N$. In that case $h^{*}$ is a polynomial which does not depend on $x_0$, and
\[
h^{*}(p)=\sum_{|\alpha|<N}\frac{\delta^\alpha(g)(x_0)}{\alpha!}\,\bigl(p-\grad h(x_0)\bigr)^{\alpha},\qquad(\grad h)^{-1}=\grad h^{*}.
\]
\end{theorem}

Theorem~\ref{thm:intro-member} combines Lemma~\ref{lem:dg} and Theorems~\ref{thm:primitive} and~\ref{thm:member}, and Theorem~\ref{thm:intro-tower} is Theorem~\ref{thm:tower}. The two theorems yield criteria of different types. Theorem~\ref{thm:intro-member} is an algebraic criterion. Condition (b) is decided by the membership algorithm of Shannon and Sweedler, which determines $h^{*}$ whenever $\grad h$ is invertible (Remark~\ref{rem:SS}). Theorem~\ref{thm:intro-tower} is a criterion of finite order in terms of Nambu derivations, and it gives the Taylor coefficients of $h^{*}$ at every point. The identity $\delta_i(g)=X_i$ also reduces the $n$ conditions of the criterion of \cite{ACH19} to the single condition $1\otimes g\in\im\varphi$ (Theorem~\ref{thm:single}). This reduces the number of conditions, not the order $N$ of the tower (Remark~\ref{rem:compare}). For $n=3$ the symmetry of the Hessian gives $\deg(\grad h)^{-1}\le2\deg h-3$ and lowers the order to $2\deg h-1$, and both bounds are attained (Proposition~\ref{prop:deg3}).

In Section~\ref{sec:pascal} we apply the inversion series of \cite{Ad26} to the identity $h^{*}\circ\grad h=g$ and write the Legendre transform as the sum of the formal Pascal series $\sum_k(-1)^k\Delta^k(g)$. Using a theorem of de Bondt \cite{dB15}, we show that in at most three variables the normal form of a gradient Keller map is $X+H$ with $JH$ strongly nilpotent, hence Pascal finite (Theorem~\ref{thm:dim3}). We give a Keller potential $\phi_0$ in four variables whose gradient is an automorphism with a normal form which is not Pascal finite, and whose inverse is determined by the membership algorithm. Section~\ref{sec:doubling} treats the doubled potential $\Phi_F(x,y)=\ip{y,F(x)}$ of a Keller map $F$, which enters the implication $\HC_{2n}\Rightarrow\JC_n$. Its formal Legendre transform is $\ip{u,G_F(v)}$, where $G_F$ is the formal inverse of $F$, and for the Nambu-tower criterion the doubling introduces no vanishing conditions beyond those for $F$ (Corollary~\ref{cor:JCsingle}). The cotangent lift of vector fields, used by the second author and Mohseni in a differential Galois study of non-integrable gradient systems \cite{HM24}, appears there as the Hamiltonian vector field of the nonlinear part of $\Phi_F$. Section~\ref{sec:scope} shows by the five-variable potential of \cite{MY26} that the properties of Keller potentials used in this paper do not imply $\HC_n$ for $n\ge5$. None of the results settles $\JC_2$ or $\HC_4$.

\section{Nambu derivations}\label{sec:nambu}

Let $F$ be a Keller map and $R=\kk[X]$. The Nambu derivations of $F$ are defined by
\begin{equation}\label{eq:nambu}
(\delta_1,\dots,\delta_n)^{\mathsf T}=((JF)^{-1})^{\mathsf T}(\partial_1,\dots,\partial_n)^{\mathsf T},
\qquad\text{that is}\qquad
\delta_i=\sum_{k=1}^n ((JF)^{-1})_{ki}\,\partial_k ,
\end{equation}
see \cite{CH11,ACH19}. Since $(JF)^{-1}=(\det JF)^{-1}\adj JF$ has polynomial entries, $\delta_i(R)\subset R$.

\begin{lemma}\label{lem:nambu}
\begin{romanlist}
\item $\delta_i(F_j)=\delta_{ij}$ for all $i,j$.
\item For every formal power series $P\in\kk[[p-p_0]]$ and every $x_0$ with $F(x_0)=p_0$ one has $\delta_i(P\circ F)=(\partial_iP)\circ F$ in $\kk[[X-x_0]]$.
\item The derivations $\delta_1,\dots,\delta_n$ commute pairwise.
\item Each $\delta_i$ has zero divergence, $\sum_k\partial_k\bigl(((JF)^{-1})_{ki}\bigr)=0$.
\end{romanlist}
\end{lemma}
\begin{proof}
(ii) By the chain rule $\partial_k(P\circ F)=\sum_j(\partial_jP\circ F)(JF)_{jk}$, hence
$\delta_i(P\circ F)=\sum_j(\partial_jP\circ F)\sum_k(JF)_{jk}((JF)^{-1})_{ki}=(\partial_iP)\circ F$. Item (i) is the case $P=p_j$.
(iii) The substitution $P\mapsto P\circ F$ is an isomorphism $\kk[[p-p_0]]\to\kk[[X-x_0]]$, since the linear part of $F$ at $x_0$ is invertible. By (ii) it transports $\partial_i$ to $\delta_i$ on $\kk[[X-x_0]]$, and the $\partial_i$ commute. As $R\subset\kk[[X-x_0]]$, the restrictions to $R$ commute. A derivation of $R$ extends uniquely to $\kk(X)$, so the extensions to $\kk(X)$ commute as well.
(iv) $(JF)^{-1}=c^{-1}\adj(JF)$ with $c=\det JF$, so it suffices to prove the Piola identity $\sum_k\partial_k\adj(JF)_{ki}=0$ for $1\le i\le n$. Let $M_{ik}$ be the minor of $JF$ obtained by deleting row $i$ and column $k$, so that $\adj(JF)_{ki}=(-1)^{i+k}M_{ik}$. Consider the $(n-1)$-form, closed as a wedge product of exact forms,
\begin{align*}
\omega_i&=dF_1\wedge\dots\wedge\widehat{dF_i}\wedge\dots\wedge dF_n\\
&=\sum_kM_{ik}\,dX_1\wedge\dots\wedge\widehat{dX_k}\wedge\dots\wedge dX_n .
\end{align*}
Then
\begin{align*}
0=d\omega_i&=\sum_k(-1)^{k-1}\partial_kM_{ik}\,dX_1\wedge\dots\wedge dX_n\\
&=(-1)^{i+1}\sum_k\partial_k\adj(JF)_{ki}\,dX_1\wedge\dots\wedge dX_n .\qedhere
\end{align*}
\end{proof}

\begin{remark}\label{rem:nambu}
The commutativity of the $\delta_i$, which is implicit in the Picard--Vessiot setting of \cite{CH11,ACH19}, is proved in (iii) without a separability argument. Statement (i) and the definition \eqref{eq:nambu} make sense over any field in which $\det JF$ is a unit. By (iv), over $\mathbb R$ the operators $-\sqrt{-1}\,\delta_i$ are formally symmetric on compactly supported smooth functions, and by (i) and (iii) the assignment $X_i\mapsto F_i$, $\partial_j\mapsto\delta_j$ defines an endomorphism of the Weyl algebra $A_n(\kk)$. This is the classical construction relating $\JC_n$ to the Dixmier conjecture \cite{BKK07,Tsu05}. The terminology ``Nambu derivations'' in \cite{CH11} refers to the generalized Hamiltonian dynamics of Nambu \cite{Nam73}. By Cramer's rule,
\[
\delta_i(P)=(\det JF)^{-1}\{F_1,\dots,F_{i-1},P,F_{i+1},\dots,F_n\},
\]
where $\{f_1,\dots,f_n\}=\det J(f_1,\dots,f_n)$ is the Nambu bracket \cite{Nam73,Tak94}. Thus $\delta_i$ is, up to a constant factor, the Nambu--Hamiltonian vector field of the $n-1$ functions $F_j$, $j\neq i$, and (iv) is the preservation of volume by such fields \cite{Nam73,Tak94}. If $F$ is only assumed to satisfy $c=\det JF\neq0$, then \eqref{eq:nambu} defines derivations of $\kk(X)$, and the Piola identity gives $\operatorname{div}\delta_i=\sum_k\partial_k\bigl(c^{-1}\adj(JF)_{ki}\bigr)=-\delta_i(c)/c$. Since $\partial_k=\sum_i(JF)_{ik}\delta_i$, all $\delta_i$ have zero divergence if and only if $c\in\kk^\times$, that is, if and only if $F$ is a Keller map. For a gradient $F=\grad h$ the derivation $\delta_i$ is the derivative with respect to $p_i=\partial_ih$ written in the coordinates $X$.
\end{remark}

\section{The criterion of Adamus, Crespo and Hajto and its converse}\label{sec:criterion}

Let $\delta_i$ act on $\kk(X)\otimes_{\kk}R$ by $\delta_i(a\otimes b)=\delta_i(a)\otimes b+a\otimes\delta_i(b)$, let $E$ be the ring of constants of $\kk(X)\otimes_\kk R$, and let
\[
\varphi\colon \kk(X)\otimes_{\kk}E\to \kk(X)\otimes_{\kk}R,\qquad \varphi(a\otimes e)=(a\otimes1)\,e .
\]
By \cite[Proposition~2.2]{ACH19} the map $\varphi$ is injective, and by \cite[Theorem~3.1]{ACH19} the surjectivity of $\varphi$ implies that $F$ has a polynomial inverse. The equivalence below is stated in \cite{ACH19}, and we include a proof of the converse direction.

\begin{lemma}\label{lem:gen}
\begin{romanlist}
\item $\im\varphi$ is a $\kk$-subalgebra of $\kk(X)\otimes R$ containing $\kk(X)\otimes1$ and $E$, and it is stable under $\delta_1,\dots,\delta_n$.
\item The elements $u_j=F_j\otimes1-1\otimes F_j$ lie in $E$.
\item $\varphi$ is surjective if and only if $1\otimes X_i\in\im\varphi$ for $1\le i\le n$.
\end{romanlist}
\end{lemma}
\begin{proof}
(i) $\varphi$ is a homomorphism of $\kk$-algebras, $\varphi(a\otimes1)=a\otimes1$ and $\varphi(1\otimes e)=e$. Moreover $\delta_i((a\otimes1)e)=(\delta_ia\otimes1)e$ because $\delta_ie=0$.
(ii) $\delta_i(u_j)=\delta_{ij}-\delta_{ij}=0$ by Lemma~\ref{lem:nambu}(i).
(iii) $\kk(X)\otimes R$ is generated as a $\kk$-algebra by $\kk(X)\otimes1$ and $1\otimes X_1,\dots,1\otimes X_n$. Apply (i).
\end{proof}

\begin{proposition}\label{prop:converse}
A Keller map $F$ has a polynomial inverse if and only if $\varphi$ is bijective.
\end{proposition}
\begin{proof}
One direction is \cite[Theorem~3.1]{ACH19}. Conversely, let $G$ be the polynomial inverse of $F$, so that $X_i=G_i(F)$. Since $1\otimes F_j=F_j\otimes1-u_j$ in $\kk(X)\otimes R$, one has, in the algebra $\kk(X)\otimes R$,
\[
1\otimes X_i=G_i(F_1\otimes1-u_1,\dots,F_n\otimes1-u_n)\in\im\varphi
\]
by Lemma~\ref{lem:gen}(i),(ii). Lemma~\ref{lem:gen}(iii) and \cite[Proposition~2.2]{ACH19} conclude.
\end{proof}

\section{Gradient maps and the Legendre potential}\label{sec:gradient}

Let $h$ be a Keller potential, $F=\grad h$, $\mathcal H=\Hess h=JF$, and
\begin{equation}\label{eq:g}
g:=\ip{X,\grad h}-h=\sum_{i=1}^nX_i\,\partial_ih-h\in\kk[X].
\end{equation}
We call $g$ the Legendre potential of $h$.

\begin{lemma}\label{lem:dg}
$\delta_i(g)=X_i$ for $1\le i\le n$.
\end{lemma}
\begin{proof}
$\partial_kg=\sum_j\mathcal H_{kj}X_j$. Since $\mathcal H$ is symmetric, $\delta_i=\sum_k(\mathcal H^{-1})_{ik}\partial_k$, and $\delta_i(g)=\sum_j(\mathcal H^{-1}\mathcal H)_{ij}X_j=X_i$.
\end{proof}

\begin{remark}\label{rem:anyring}
The computation uses only the symmetry of $\Hess h$ and no division by integers. Lemma~\ref{lem:dg} holds over any commutative ring in which $\det\Hess h$ is a unit. The remaining results use characteristic zero at the following points. The criterion of \cite{ACH19}, and hence Theorem~\ref{thm:single}, rests on Picard--Vessiot theory in characteristic zero. The proof of Theorem~\ref{thm:primitive} uses $P'(g)\neq0$, which fails when $g$ is inseparable over $\kk(\grad h)$. Formula \eqref{eq:explicit} divides by $\alpha!$, and the step from the polynomiality of $\grad h^{*}$ to that of $h^{*}$ in the proof of Theorem~\ref{thm:tower} fails in characteristic $p$, since the series $\sum_{k\ge1}X_1^{pk}$ has gradient zero.
\end{remark}

\begin{theorem}[single condition]\label{thm:single}
$\grad h$ is a polynomial automorphism if and only if $1\otimes g\in\im\varphi$.
\end{theorem}
\begin{proof}
If $1\otimes g\in\im\varphi$, then $1\otimes X_i=\delta_i(1\otimes g)\in\im\varphi$ by Lemmas~\ref{lem:gen}(i) and~\ref{lem:dg}, so $\varphi$ is surjective. Proposition~\ref{prop:converse} gives both directions.
\end{proof}

\subsection{The Legendre transform}
Fix $x_0\in\kk^n$ and put $p_0=\grad h(x_0)$. Since $\mathcal H(x_0)$ is invertible, $\grad h$ has a unique formal inverse $G\in\kk[[p-p_0]]^n$ with $G(p_0)=x_0$. The formal Legendre transform of $h$ at $x_0$ is
\[
h^{*}(p):=\ip{p,G(p)}-h(G(p))\in\kk[[p-p_0]] .
\]
By definition $h^{*}\circ\grad h=g$. The series $h^{*}$ depends on $x_0$. If $\grad h$ is a polynomial automorphism, then $G$ is the Taylor expansion at $p_0$ of the polynomial inverse, so $h^{*}=\ip{p,(\grad h)^{-1}(p)}-h\bigl((\grad h)^{-1}(p)\bigr)$ is the same polynomial for every $x_0$. Statements on $h^{*}$ without mention of a base point hold for every fixed $x_0$.

\begin{lemma}[{\cite{Men06}}]\label{lem:leg}
$\grad h^{*}=G$.
\end{lemma}
\begin{proof}
$\partial_{p_i}h^{*}=G_i+\sum_j\bigl(p_j-(\partial_jh)(G)\bigr)\partial_{p_i}G_j=G_i$, because $(\grad h)(G(p))=p$.
\end{proof}

\begin{remark}[relation to the work of Meng]\label{rem:meng}
The Hessian conjecture was formulated by Meng through the formal Legendre transform \cite{Men06}, and independently, without the Legendre transform, by de Bondt and van den Essen \cite{dBvdE05}, see the acknowledgement in \cite{Men06}. Meng proves that $\HC_n$ for all $n$ is equivalent to $\JC_n$ for all $n$ \cite[Proposition~1.4]{Men06}. In a fixed dimension his argument gives $\HC_{2n}\Rightarrow\JC_n\Rightarrow\HC_n$, and in particular it does not give $\JC_2\Rightarrow\HC_4$. He also proves $\HC_n$ over $\mathbb R$ when $\Hess h$ is definite at one point \cite[Theorem~1.3]{Men06}, and gives a tree formula for the Taylor coefficients of $h^{*}$ at $\grad h(0)$ in terms of those of $h$ at $0$ \cite[Theorem~2.1]{Men06}. Proposition~\ref{prop:tower} below expresses the Taylor coefficients of $h^{*}$ at $\grad h(x_0)$, for every $x_0$, as values at $x_0$ of the polynomials $\delta^\alpha(g)$.
\end{remark}

\begin{proposition}\label{prop:tower}
For every multi-index $\alpha$, the identity $\delta^\alpha(g)=(\partial^\alpha h^{*})\circ\grad h$ holds in $\kk[[X-x_0]]$. In particular, $(\partial^\alpha h^{*})\circ\grad h$ is a polynomial.
\end{proposition}
\begin{proof}
Induction on $|\alpha|$. The case $\alpha=0$ is $h^{*}\circ\grad h=g$. The inductive step is Lemma~\ref{lem:nambu}(ii) applied to $P=\partial^\alpha h^{*}$, the order of the derivations being irrelevant by Lemma~\ref{lem:nambu}(iii).
\end{proof}

\begin{theorem}[Nambu tower of the Legendre potential]\label{thm:tower}
Let $d=\deg h\ge2$ and $N=(d-1)^{n-1}+2$. The following are equivalent.
\begin{alphalist}
\item $\grad h$ is a polynomial automorphism.
\item $h^{*}$ is a polynomial for some, equivalently for every, $x_0\in\kk^n$.
\item $\delta^\alpha(g)=0$ for all $\alpha$ with $|\alpha|=N$.
\item $D^N(g)=0$, where $D=\sum_iY_i\delta_i$ on $\kk[X,Y]$ with $\delta_i(Y_j)=0$.
\end{alphalist}
Moreover, if $\delta^\alpha(g)=0$ for all $|\alpha|=N'$ for some $N'\ge2$, then (a)--(d) hold. For $N'=1$ this hypothesis never holds, since $\delta_i(g)=X_i$. The order $N$ comes from the bound on the degree of the inverse. For a given $h$ the tower may vanish at a lower order. If (a)--(d) hold, then for every $x_0\in\kk^n$
\begin{equation}\label{eq:explicit}
h^{*}(p)=\sum_{|\alpha|<N}\frac{\delta^\alpha(g)(x_0)}{\alpha!}\bigl(p-\grad h(x_0)\bigr)^\alpha,\qquad (\grad h)^{-1}=\grad h^{*}.
\end{equation}
\end{theorem}
\begin{proof}
We prove (a)$\Rightarrow$(b) with the bound $\deg h^{*}<N$, which uses Lemma~\ref{lem:leg} and the bound on the degree of the inverse, then (b)$\Rightarrow$(a) by Lemma~\ref{lem:leg}, then (a)$\Rightarrow$(c) and (c)$\Rightarrow$(b) by Proposition~\ref{prop:tower}, and finally (c)$\Leftrightarrow$(d), which is formal.
(a)$\Rightarrow$(b). The inverse coincides with $G$. Since $\deg\grad h=d-1$, one has $\deg G\le(d-1)^{n-1}$ by \cite[Chapter~I, Corollary~1.4]{BCW82}. By Lemma~\ref{lem:leg}, the formal power series $h^{*}$ has the polynomial gradient $G$, so $h^{*}$ is a polynomial, and $\deg h^{*}\le\deg G+1<N$.
(b)$\Rightarrow$(a). By Lemma~\ref{lem:leg}, $G=\grad h^{*}$ is a polynomial map. The identities $G\circ\grad h=X$ and $\grad h\circ G=p$ hold in $\kk[[X-x_0]]$ and in $\kk[[p-p_0]]$. Both sides of each identity are polynomials, and a polynomial which vanishes as a formal power series at a point is zero, so they are polynomial identities.
(a)$\Rightarrow$(c). By the proof of (a)$\Rightarrow$(b), $\deg h^{*}<N$, so $\partial^\alpha h^{*}=0$ for $|\alpha|=N$. Apply Proposition~\ref{prop:tower}.
(c)$\Rightarrow$(b). By Proposition~\ref{prop:tower}, $(\partial^\alpha h^{*})\circ\grad h=0$ for $|\alpha|=N$. Since the linear part $\mathcal H(x_0)$ of $\grad h$ at $x_0$ is invertible, $P\mapsto P\circ\grad h$ is an isomorphism $\kk[[p-p_0]]\to\kk[[X-x_0]]$, as in the proof of Lemma~\ref{lem:nambu}(iii), so $\partial^\alpha h^{*}=0$ for $|\alpha|=N$ and $h^{*}$ is a polynomial of degree $<N$.
The proof of (c)$\Rightarrow$(b) uses only the vanishing at one order, which gives the last assertion.
(c)$\Leftrightarrow$(d). By Lemma~\ref{lem:nambu}(iii), $D^N(g)=\sum_{|\alpha|=N}\frac{N!}{\alpha!}Y^\alpha\delta^\alpha(g)$, and the monomials $Y^\alpha$ are linearly independent over $\kk[X]$.
Formula \eqref{eq:explicit} is the Taylor expansion of the polynomial $h^{*}$ at $p_0$, whose coefficients are $(\partial^\alpha h^{*})(p_0)=\delta^\alpha(g)(x_0)$ by Proposition~\ref{prop:tower}.
\end{proof}

\begin{remark}\label{rem:compare}
\begin{arabiclist}
\item Since $X_i=\delta_i(g)$, the condition $D^N(g)=0$ implies $D^N(X_i)=\delta_i(D^N(g))=0$, so $D$ is locally nilpotent. Indeed, by the Leibniz formula $D^m(PQ)=\sum_{k}\binom mkD^k(P)D^{m-k}(Q)$ the elements of $\kk[X,Y]$ annihilated by a power of $D$ form a subalgebra, and it contains $X_1,\dots,X_n$ and $Y_1,\dots,Y_n$. In the gradient case the criterion of \cite[Proposition~3.1.4]{vdE00}, see also \cite[Remark~3.3]{ACH19}, needs to be checked on the single element $g$, one order higher than on the coordinates.
\item The reduction concerns the number of conditions, which drops from $n$ to one, and not the cost of the computation. It replaces $n$ targets of degree one by one target of degree $\deg h$, it does not change the number $\binom{N+n-1}{n-1}$ of multi-indices of order $N$, and it does not remove the computation of the $\delta_i$ and of the powers of $D$, observed in \cite[Remark~3.3]{ACH19} to dominate the running time. In the membership test of Remark~\ref{rem:SS} below, the Gr\"obner basis is the same for $g$ as for the coordinates, and the reduction concerns the number of normal form computations.
\item The coefficient of $(p-\grad h(x_0))^\alpha$ in $h^{*}$ is $\delta^\alpha(g)(x_0)/\alpha!$. Formula \eqref{eq:explicit} therefore determines the inverse from finitely many values at a single point.
\end{arabiclist}
\end{remark}

The following lemma is a standard consequence of the uniqueness of the formal inverse.

\begin{lemma}[window of zeros]\label{lem:window}
Let $F=X+H\in\kk[X]^n$ with $\operatorname{ord}H\ge2$, $e=\deg F\ge2$, and let $G=\sum_{s\ge1}G_s$ be its formal inverse. Let $B\ge1$. If $G_s=0$ for $B<s\le eB$, then $G=G_{\le B}$ is a polynomial and $F$ is an automorphism.
\end{lemma}
\begin{proof}
Put $R=G-G_{\le B}$, so $\operatorname{ord}R>eB$. From $G\circ F=X$ we get $R\circ F=X-G_{\le B}\circ F$. Since $F=X+H$ with $\operatorname{ord}H\ge2$, composition with $F$ preserves the order of a nonzero formal power series, the lowest homogeneous part of $Q\circ F$ being that of $Q$. The right-hand side is a polynomial of degree $\le eB$, and the left-hand side has order $>eB$. Hence $R\circ F=0$ and $R=(R\circ F)\circ G=0$.
\end{proof}

\begin{corollary}[pointwise test]\label{cor:pointwise}
Let $h$ be a Keller potential of degree $d\ge3$, $x_0\in\kk^n$, $p_0=\grad h(x_0)$ and $B=(d-1)^{n-1}$. Then $\grad h$ is an automorphism if and only if the homogeneous parts of $h^{*}$ in $p-p_0$ of degrees $B+2,\dots,(d-1)B+1$ vanish, that is, if and only if $\delta^\alpha(g)(x_0)=0$ for all $\alpha$ with $B+2\le|\alpha|\le(d-1)B+1$.
\end{corollary}
\begin{proof}
By Proposition~\ref{prop:tower}, $\delta^\alpha(g)(x_0)=(\partial^\alpha h^{*})(p_0)$, which gives the equivalence of the two formulations. By Lemma~\ref{lem:leg} the part $h^{*}_{s+1}$ of degree $s+1$ and the part $G_s$ of degree $s$ of $G-x_0$ determine each other by $G_s=\grad h^{*}_{s+1}$ and $(s+1)h^{*}_{s+1}=\ip{p-p_0,G_s}$. Hence the hypothesis says that $G_s=0$ for $B<s\le(d-1)B$. Put $\tilde F(Y)=\mathcal H(x_0)^{-1}\bigl(\grad h(x_0+Y)-p_0\bigr)$. Then $\tilde F=Y+\tilde H$ with $\operatorname{ord}\tilde H\ge2$, and $\deg\tilde F=d-1$, since the top homogeneous part of $h$ has a nonzero gradient. From $G\circ\grad h=X$ and $\grad h\circ G=p$ one checks that $\tilde G(q)=G(p_0+\mathcal H(x_0)q)-x_0$ satisfies $\tilde G\circ\tilde F=Y$ and $\tilde F\circ\tilde G=q$, so it is the formal inverse of $\tilde F$, and the homogeneous parts are $\tilde G_s(q)=G_s(\mathcal H(x_0)q)$. Hence $\tilde G_s=0$ for $B<s\le(d-1)B$, and Lemma~\ref{lem:window} applied to $\tilde F$ shows that $\tilde F$, and therefore $\grad h$, is an automorphism. The converse follows from $\deg G\le B$ \cite{BCW82}.
\end{proof}

The test is pointwise in the sense that only values at the single point $x_0$ are used. It does not reduce the number of multi-indices to be considered.

\begin{remark}\label{rem:d2}
For $d=2$ the map $\grad h$ is affine with invertible linear part, so the assumption $d\ge3$ is no restriction.
\end{remark}

\begin{proposition}\label{prop:third}
\begin{romanlist}
\item $\delta_i\delta_j(g)=(\mathcal H^{-1})_{ij}$.
\item $\delta_k\bigl((\mathcal H^{-1})_{ij}\bigr)=\bigl(\partial_{p_i}\partial_{p_j}\partial_{p_k}h^{*}\bigr)\circ\grad h$.
\item $\delta_k(\mathcal H_{ij})=\sum_l(\mathcal H^{-1})_{kl}\,\partial_i\partial_j\partial_lh$.
\end{romanlist}
\end{proposition}
\begin{proof}
(i) $\delta_i(X_j)=(\mathcal H^{-1})_{ij}$ by \eqref{eq:nambu}. (ii) Combine (i) with Proposition~\ref{prop:tower}. (iii) Definition \eqref{eq:nambu}.
\end{proof}

\begin{remark}\label{rem:hessian-structures}
In the language of Hessian structures \cite{Shi07}, $h$ and $h^{*}$ are dual potentials, the Hessian metrics are $\mathcal H$ and $\Hess h^{*}=\mathcal H^{-1}\circ G$, and $T_{ijl}=\partial_i\partial_j\partial_lh$ is the cubic form. Proposition~\ref{prop:third}(iii) expresses the Nambu derivative of the metric through $T$ with one index raised by $\mathcal H^{-1}$, and (ii) expresses the Nambu derivative of the dual metric through the cubic form of $h^{*}$.
\end{remark}

\begin{example}\label{ex:n3}
For $h=\frac16X_3^3+X_1X_3+\frac12X_2^2$ one has $\det\Hess h=-1$, $h^{*}=-\frac16p_1^3+p_1p_3+\frac12p_2^2$, $D^4(g)=0$ and $D^3(g)\neq0$, while $D^3(X_1)=0$ and $D^2(X_1)\neq0$ and $D^2(X_2)=D^2(X_3)=0$.
For $h=\frac16X_4^3+X_3X_4+X_1X_2+X_1^6$ one has $\det\Hess h=1$, $D^7(g)=0$ and $D^6(g)\neq0$, and \eqref{eq:explicit} at $x_0=(1,-2,3,1)$ gives $h^{*}=p_1p_2-p_2^6+p_3p_4-\frac16p_3^3$, so that $(\grad h)^{-1}=\grad h^{*}=\bigl(p_2,\,p_1-6p_2^5,\,p_4-\frac12p_3^2,\,p_3\bigr)$.
\end{example}

\subsection{The Legendre potential as a primitive element}

In this subsection $h\in\kk[X]$ is only assumed to have a Hessian determinant which is a nonzero polynomial, not necessarily constant, and $F=\grad h$. The derivations $\delta_i$ of \eqref{eq:nambu} are then derivations of $\kk(X)$, and the proofs of Lemma~\ref{lem:nambu}(i)--(iii) and of Lemma~\ref{lem:dg} apply verbatim at any point $x_0$ with $\det\Hess h(x_0)\neq0$. Since $\det JF\neq0$, the polynomials $F_1,\dots,F_n$ are algebraically independent over $\kk$ by the Jacobian criterion. Hence $\kk(F)$ has transcendence degree $n$ over $\kk$, the field $\kk(X)$ is algebraic over $\kk(F)$, and since $\kk(X)$ is finitely generated, the extension $\kk(X)\supset\kk(F)$ is finite. Theorem~\ref{thm:primitive} below shows that $g$ generates it. We put $\mu(F)=[\kk(X):\kk(F)]$, the number of points of a general fibre of $F$.

\begin{theorem}[primitive element]\label{thm:primitive}
Let $h\in\kk[X]$ such that $\det\Hess h$ is a nonzero polynomial, $F=\grad h$ and $g=\ip{X,\grad h}-h$. Then $\kk(X)=\kk(F)(g)$. In particular, the minimal polynomial of $g$ over $\kk(F)$ has degree $\mu(F)$.
\end{theorem}
\begin{proof}
For $Q\in\kk(T_1,\dots,T_n)$ one has $\delta_i(Q(F))=(\partial_{T_i}Q)(F)$ by the chain rule and Lemma~\ref{lem:nambu}(i), so $\kk(F)$ is stable under the $\delta_i$, as observed in \cite{CH11}. Let $P\in\kk(F)[T]$ be the minimal polynomial of $g$ over $\kk(F)$, and for $P=\sum_ja_jT^j$ with $a_j\in\kk(F)$ put $P^{\delta_i}=\sum_j\delta_i(a_j)T^j$. Then $0=\delta_i(P(g))=P^{\delta_i}(g)+P'(g)\,\delta_i(g)$, and $P'(g)\neq0$, since $P$ is irreducible and, in characteristic zero, prime to its nonzero derivative $P'$ of smaller degree. Hence $X_i=\delta_i(g)=-P^{\delta_i}(g)/P'(g)\in\kk(F)(g)$ by Lemma~\ref{lem:dg}.
\end{proof}

\begin{remark}\label{rem:primitive}
Theorem~\ref{thm:primitive} identifies the Legendre potential as a primitive element of $\kk(X)\supset\kk(\grad h)$. For a Keller potential it gives the following reading of known criteria. By \cite[Theorem~2]{CH11} and Campbell's theorem \cite{Cam73}, $\grad h$ is an automorphism if and only if $\kk(X)\supset\kk(\grad h)$ is Galois, that is, if and only if the minimal polynomial of $g$ over $\kk(\grad h)$ splits in $\kk(X)$. By Keller's theorem, $\grad h$ is an automorphism if and only if this polynomial has degree one. For $h=X_1^3+X_2^3$, which is not a Keller potential, the minimal polynomial of $g=2X_1^3+2X_2^3$ over $\kk(\grad h)$ has degree $4=\mu(\grad h)$, since its conjugates $2(\pm X_1^3\pm X_2^3)$ are pairwise distinct.
\end{remark}

\subsection{Algebra membership}

\begin{theorem}[membership]\label{thm:member}
Let $h$ be a Keller potential, $F=\grad h$ and $g=\ip{X,\grad h}-h$. The following are equivalent.
\begin{alphalist}
\item $\grad h$ is a polynomial automorphism.
\item $g\in\kk[\partial_1h,\dots,\partial_nh]$.
\item $g\in\kk(\partial_1h,\dots,\partial_nh)$.
\end{alphalist}
If these hold and $g=P(\partial_1h,\dots,\partial_nh)$ with $P\in\kk[T_1,\dots,T_n]$, then $P=h^{*}$ for every base point $x_0$.
\end{theorem}
\begin{proof}
(a)$\Rightarrow$(b). $\kk[F]=\kk[X]$. (b)$\Rightarrow$(c) is clear.
(c)$\Rightarrow$(a). By Theorem~\ref{thm:primitive}, $\kk(X)=\kk(F)(g)=\kk(F)$, and a Keller map with this property is a polynomial automorphism \cite{Kel39}, \cite[Chapter~I, Theorem~2.1]{BCW82}.
Finally, $P$ is unique, since $F_1,\dots,F_n$ are algebraically independent. If $g=P(F)$, then $P\circ F=h^{*}\circ F$ in $\kk[[X-x_0]]$, and composition with $F$ is injective, so $P=h^{*}$.
\end{proof}

\begin{remark}\label{rem:SS}
Condition (b) is decided by the algorithm of \cite{SS88}. Let $\mathcal G$ be a Gr\"obner basis of the ideal $(\partial_1h-T_1,\dots,\partial_nh-T_n)\subset\kk[X,T]$ for a block order with $X\gg T$, and let $r$ be the remainder of $g$ on division by $\mathcal G$. Since $\partial_1h,\dots,\partial_nh$ are algebraically independent, the ideal meets $\kk[T]$ only in $0$, and by \cite{SS88} condition (b) holds if and only if $r\in\kk[T]$, in which case $g=r(\partial_1h,\dots,\partial_nh)$ and $r=h^{*}$ by Theorem~\ref{thm:member}. Indeed, if $g=P(\partial_1h,\dots,\partial_nh)$, then $g-P(T)$ lies in the ideal, $P(T)$ is reduced for the block order, and the uniqueness of remainders gives $r=P(T)$. One Gr\"obner basis and one reduction therefore decide the invertibility of $\grad h$ and, if $\grad h$ is invertible, determine its Legendre transform. This is a decision procedure, and no bound on the cost of the Gr\"obner basis computation is implied. The Wronskian criterion of \cite{ABH17} requires the membership of $\frac12n^2(n+1)-n$ elements. For $h=X_1X_3+\frac12X_2^2+X_2X_3^2+X_3^4$, with $\det\Hess h=-1$, the remainder is $r=-\frac12T_1^4-T_1^2T_2+T_1T_3+\frac12T_2^2$.
\end{remark}

\section{Doubling and the cotangent lift}\label{sec:doubling}

For a polynomial map $F\colon\kk^n\to\kk^n$ put $\Phi_F(x,y)=\ip{y,F(x)}\in\kk[x,y]$. Then
\begin{equation}\label{eq:doubling}
\begin{gathered}
\grad\Phi_F=\bigl(JF(x)^{\mathsf T}y,\,F(x)\bigr),\qquad
\Hess\Phi_F=\begin{pmatrix}\sum_iy_i\Hess F_i & JF^{\mathsf T}\\ JF&0\end{pmatrix},\\
\det\Hess\Phi_F=(-1)^n(\det JF)^2 .
\end{gathered}
\end{equation}

For a polynomial map $H$ we put $\mathcal E(H)=JH\cdot X$, so that $\mathcal E(H)=mH$ if $H$ is homogeneous of degree $m$.

\begin{proposition}\label{prop:doubling}
Let $F$ be a Keller map and $h_0\in\kk[x]$. Put $\Phi=\Phi_F+h_0$.
\begin{romanlist}
\item $\Phi$ is a Keller potential, and $\grad\Phi$ is injective if and only if $F$ is injective.
\item If $F$ has a polynomial inverse $G$, then
\[
(\grad\Phi)^{-1}(u,v)=\Bigl(G(v),\,\bigl(JF(G(v))^{\mathsf T}\bigr)^{-1}\bigl(u-\grad h_0(G(v))\bigr)\Bigr).
\]
\end{romanlist}
In particular, $\HC_{2n}\Rightarrow\JC_n$, and $\JC_n$ is equivalent to $\HC_{2n}$ restricted to potentials affine in $y$.
\end{proposition}
\begin{proof}
The Hessian of $\Phi$ is that of \eqref{eq:doubling} with $\Hess h_0$ added to the upper left block, and its determinant is computed by the Laplace expansion along the last $n$ rows, which gives $(-1)^n(\det JF)^2$. Further $\grad\Phi=\bigl(JF^{\mathsf T}y+\grad h_0,\,F\bigr)$.
(i) Let $F(a)=F(b)$ with $a\neq b$. Since $\det JF\in\kk^\times$, the matrix $(JF^{\mathsf T})^{-1}$ has polynomial entries. Put $y'=(JF(b)^{\mathsf T})^{-1}\bigl(\grad h_0(a)-\grad h_0(b)\bigr)$. Then $\grad\Phi(a,0)=\grad\Phi(b,y')$. Conversely, if $F$ is injective and $\grad\Phi(a,y)=\grad\Phi(a',y')$, then $F(a)=F(a')$, so $a=a'$, and $JF(a)^{\mathsf T}(y-y')=0$ gives $y=y'$.
(ii) A direct substitution, using that $(JF^{\mathsf T})^{-1}=(\det JF)^{-1}\adj(JF)^{\mathsf T}$ is polynomial.
For the last assertions, let $F$ be a Keller map. Then $\Phi_F$ is a Keller potential, and if $\HC_{2n}$ holds, $\grad\Phi_F$ is an automorphism by Theorem~\ref{thm:tower}, so $F$ is injective by (i) and has a polynomial inverse \cite{BBR62}, \cite[Chapter~4]{vdE00}. Conversely, a polynomial in $(x,y)$ which is affine in $y$ has the form $h_0(x)+\ip{y,F(x)}$, and its Hessian determinant $(-1)^n(\det JF)^2$ is a nonzero constant only if $F$ is a Keller map. If $\JC_n$ holds, $F$ has a polynomial inverse, so $\grad\Phi$ is an automorphism by (ii), and $\Phi^{*}$ is a polynomial by Theorem~\ref{thm:tower}.
\end{proof}

For $F=X+H$ the symmetric lift of de Bondt and van den Essen \cite{dBvdE05}, in the form recalled in \cite{Pre26}, is $\mathrm{id}+\grad f$ with $f(x,y)=-\sqrt{-1}\,\ip{y,H(x+\sqrt{-1}\,y)}$. With $A(x,y)=(x+\sqrt{-1}\,y,-\sqrt{-1}\,y)$ one has $\frac12\ip{x,x}+\frac12\ip{y,y}+f=(\Phi_F+\frac12\ip{x,x})\circ A$, so the symmetric lift equals $A^{\mathsf T}\circ\grad(\Phi_F+\frac12\ip{x,x})\circ A$.

The cotangent lift $T^{*}F(x,u)=\bigl(F(x),(JF(x)^{\mathsf T})^{-1}u\bigr)$ satisfies $T^{*}F(x,u)=(v,y)$ if and only if $\grad\Phi_F(x,y)=(u,v)$, so that $\Phi_F$ is a generating function of $T^{*}F$. We write the coordinates of $\grad\Phi_F$ as $(u,v)=(JF^{\mathsf T}y,F)$.

\begin{proposition}\label{prop:doubled-tower}
Let $F$ be a Keller map, $\Phi=\Phi_F$, and let $g_\Phi=\ip{y,JF\,x}$ be the Legendre potential of $\Phi$. Let $\delta_{u_i},\delta_{v_j}$ be the Nambu derivations of $\grad\Phi$ and $\delta^F_j$ those of $F$.
\begin{romanlist}
\item The formal Legendre transform of $\Phi$ is $\Phi^{*}(u,v)=\ip{u,G_F(v)}$, where $G_F$ is the formal inverse of $F$.
\item For every multi-index $\beta$ in $v$ one has $\delta_{u_i}\delta_v^\beta(g_\Phi)=(\delta^F)^\beta(x_i)$ and $\delta_v^\beta(g_\Phi)=\ip{JF^{\mathsf T}y,(\delta^F)^\beta(x)}$, and $\delta^\alpha(g_\Phi)=0$ whenever $\alpha$ contains at least two derivations $\delta_{u_i}$.
\end{romanlist}
\end{proposition}
\begin{proof}
(i) The formal inverse of $\grad\Phi$ is $(u,v)\mapsto\bigl(G_F(v),(JF^{\mathsf T})^{-1}(G_F(v))\,u\bigr)$. With $x=G_F(v)$ and $y=(JF(x)^{\mathsf T})^{-1}u$ one gets $\ip{u,x}+\ip{v,y}-\ip{y,F(x)}=\ip{u,x}$.
(ii) By Proposition~\ref{prop:tower}, $\delta^\alpha(g_\Phi)=(\partial^\alpha\Phi^{*})\circ\grad\Phi$. By (i), $\partial_{u_i}\partial_v^\beta\Phi^{*}=\partial_v^\beta(G_F)_i$, $\partial_v^\beta\Phi^{*}=\ip{u,\partial_v^\beta G_F}$, and derivatives of order $\ge2$ in $u$ vanish. Finally $(\partial_v^\beta G_F)\circ F=(\delta^F)^\beta(x)$ by Lemma~\ref{lem:nambu}(ii) applied to $F$.
\end{proof}

\begin{remark}\label{rem:doubled-order}
Theorem~\ref{thm:tower} applied to $\Phi$, of degree $\deg F+1$ in $2n$ variables, gives the order $(\deg F)^{2n-1}+2$. Proposition~\ref{prop:doubled-tower}(i) and \cite{BCW82} show that $\deg\Phi^{*}\le(\deg F)^{n-1}+1$ when $F$ is invertible, so the order $(\deg F)^{n-1}+2$ suffices.
\end{remark}

\begin{theorem}[Legendre potential of the doubling]\label{thm:doublingpot}
The Legendre potential of $\Phi_F$ is $\ip{y,JF(x)\,x}$. If $F=X+H$, it equals $\ip{y,x}+\ip{y,\mathcal E(H)}$, and if moreover $H$ is homogeneous of degree $m$, it equals $\ip{y,x}+m\ip{y,H}$.
\end{theorem}
\begin{proof}
By \eqref{eq:doubling}, $\ip{x,JF^{\mathsf T}y}+\ip{y,F}-\ip{y,F}=\ip{y,JF\,x}$. Further $J(X+H)\,x=x+\mathcal E(H)$, and $\mathcal E(H)=mH$ for $H$ homogeneous of degree $m$.
\end{proof}

\begin{corollary}\label{cor:JCsingle}
Let $F$ be a Keller map, $\Phi=\Phi_F$ and $N\ge1$. Then $\delta^\alpha(g_\Phi)=0$ for all $|\alpha|=N$ if and only if $(\delta^F)^\beta(x_i)=0$ for all $i$ and all $|\beta|=N-1$. In particular, for the cubic homogeneous reduction of \cite{BCW82,Yag80}, $F=X+H$ with $H$ cubic homogeneous, the Legendre potential of $\Phi_F$ is $\ip{y,x}+3\ip{y,H}$, and Theorem~\ref{thm:tower} applied to $\Phi_F$ is the criterion of \cite[Proposition~3.1.4]{vdE00} applied to $F$, with the order raised by one.
\end{corollary}
\begin{proof}
By Proposition~\ref{prop:doubled-tower}(ii), the elements $\delta^\alpha(g_\Phi)$ with $|\alpha|=N$ are $0$, $(\delta^F)^\beta(x_i)$ with $|\beta|=N-1$, and $\ip{JF^{\mathsf T}y,(\delta^F)^\beta(x)}$ with $|\beta|=N$, according as $\alpha$ contains at least two, exactly one or no derivation $\delta_{u_i}$. The elements of the third kind vanish if those of the second kind vanish, since $(\delta^F)^\beta=\delta^F_j(\delta^F)^{\beta'}$ for some $j$ and $|\beta'|=N-1$.
\end{proof}

\begin{remark}\label{rem:symlift}
The symmetric lift was introduced in \cite{dBvdE05} to reduce $\JC_n$ to Keller maps with symmetric Jacobian matrix. Corollary~\ref{cor:JCsingle} shows that, for the criterion of Theorem~\ref{thm:tower}, the doubled potential gives no condition beyond the criterion for $F$ itself. In particular, in the setting of Nambu derivations the passage from $\JC_n$ to $\HC_{2n}$ does not simplify the verification of invertibility.
\end{remark}

\begin{remark}[the cotangent lift]\label{rem:AZ}
Let $V$ be a polynomial vector field on $\kk^n$ and $f_V(q,p)=\ip{p,V(q)}$ its momentum function on $T^{*}\kk^n$ with the form $\sum dq_i\wedge dp_i$. The Hamiltonian vector field of $f_V$ is the cotangent lift $\hat V=(V(q),-JV(q)^{\mathsf T}p)$ of \cite{AZ10}, used in \cite{HM24}, and by \cite{AZ10} the differential Galois groups of the variational equations of $V$ and of $\hat V$ are isomorphic. As in \eqref{eq:doubling}, $\det\Hess f_V=(-1)^n(\det JV)^2$, so $f_V$ is a Keller potential if and only if $V$ is a Keller map. For $F=X+H$ one has $\Phi_F=\ip{y,x}+f_H$, and by Theorem~\ref{thm:doublingpot} the Legendre potential of $\Phi_F$ is $\ip{y,x}+f_{\mathcal E(H)}$. Thus the nonlinear part of the doubled potential of a Keller map $X+H$ is the momentum function of $H$. The example of \cite[Example~3.2]{HM24} is a gradient field $\grad P$ with nonconstant $\det\Hess P$, hence not a Keller map.
\end{remark}

\section{The Pascal--Legendre series}\label{sec:pascal}

Throughout this section $F=X+H\in\kk[X]^n$ with $t:=\min_i\operatorname{ord}H_i\ge2$, and $\Delta_F(Q)=Q\circ F-Q$ for $Q\in\kk[[X]]$. Recall \cite{ABCH19} that $F$ is Pascal finite if $\Delta_F^m(X)=0$ for some $m$, and that in characteristic zero this holds if and only if $F$ is the exponential of a locally nilpotent derivation.

\begin{lemma}[{scalar form of \cite[Lemma~2, Theorem~5]{Ad26}}]\label{lem:scalar}
For $Q\in\kk[[X]]$ with $Q(0)=0$ one has $\operatorname{ord}\Delta_F^kQ\ge\operatorname{ord}Q+k(t-1)$. The series $\sum_{k\ge0}(-1)^k\Delta_F^kQ$ converges in the $(X)$-adic topology and equals $Q\circ F^{-1}$, where $F^{-1}$ is the formal inverse.
\end{lemma}
\begin{proof}
For a monomial $M$ of degree $r$ every term of $M\circ F-M$ has degree $\ge r+t-1$. Applying this to the terms of $Q$ shows that every term of $\Delta_F(Q)$ has degree at least $\operatorname{ord}Q+t-1$. Hence $\Delta_F$ raises the order by at least $t-1\ge1$, the operator $\sum_k(-\Delta_F)^k$ converges and is inverse to $\sigma_F=\mathrm{id}+\Delta_F\colon Q\mapsto Q\circ F$. On the other hand $\sigma_F^{-1}=\sigma_{F^{-1}}$.
\end{proof}

The series of Lemma~\ref{lem:scalar} always converges in the $(X)$-adic topology. Its termination is an additional property. Applied to the Legendre potential, termination is sufficient for the invertibility of $\grad h$ (Corollary~\ref{cor:single-pascal}) but not necessary (Proposition~\ref{prop:phi0-g}).

\begin{proposition}[Pascal--Legendre series]\label{prop:PL}
Let $h$ be a Keller potential in normal form, that is, $h(0)=0$, $\grad h(0)=0$ and $\Hess h(0)=I$, so that $\grad h=X+H$ with $t\ge2$ and $g(0)=0$. Then
\[
h^{*}=\sum_{k\ge0}(-1)^k\Delta_{\grad h}^k(g),\qquad g=\ip{X,\grad h}-h=\tfrac12\ip{X,X}+\sum_{j\ge3}(j-1)h_j,
\]
where $h=\frac12\ip{X,X}+\sum_{j\ge3}h_j$ is the decomposition into homogeneous parts.
\end{proposition}
\begin{proof}
$h^{*}\circ\grad h=g$, so $h^{*}=g\circ G$. Apply Lemma~\ref{lem:scalar}. The formula for $g$ follows from $\ip{X,\grad h_j}=jh_j$.
\end{proof}

\begin{remark}\label{rem:zhao}
For $h=\frac12\ip{X,X}-P$ with $P$ homogeneous of degree $4$ and $\Hess P$ nilpotent, Zhao \cite{Zha05,Zha07} writes the inverse of $X-\grad P$ as $X+\grad Q$ with $Q=\sum_{m\ge0}\frac{\triangle^m(P^{m+1})}{2^m m!(m+1)!}$, and the vanishing conjecture $\triangle^m(P^{m+1})=0$ for $m\gg0$ is the statement that this series terminates. Proposition~\ref{prop:PL} is a difference analogue valid for every Keller potential. After the counterexamples to $\JC_n$, an explicit Hessian nilpotent quartic in $38$ variables whose gradient map $X-\grad P$ is not injective is given in \cite{Pre26}. For such $P$ neither Zhao's series nor the series of Proposition~\ref{prop:PL} terminates.
\end{remark}

\begin{remark}\label{rem:normalform}
Let $h$ be a Keller potential, $x_0\in\kk^n$, $L=\Hess h(x_0)$, and $A\in\mathrm{GL}_n(\kk)$ with $A^{\mathsf T}LA=I$. Then $\tilde h(X)=h(x_0+AX)-h(x_0)-\ip{\grad h(x_0),AX}$ is in normal form, and $\grad\tilde h=A^{-1}\circ\tilde F\circ A$ with $\tilde F(Y)=L^{-1}\bigl(\grad h(x_0+Y)-\grad h(x_0)\bigr)$, since $A^{\mathsf T}=A^{-1}L^{-1}$. We call $\grad\tilde h$ a normal form of $\grad h$ at $x_0$. Pascal finiteness, the Pascal depth and the degrees of the iterates and of the inverse are invariant under linear conjugation.
\end{remark}

\begin{corollary}\label{cor:single-pascal}
Let $h$ be a Keller potential in normal form. If $\Delta_{\grad h}^m(g)=0$ for some $m$, then $h^{*}=\sum_{k<m}(-1)^k\Delta^k_{\grad h}(g)$ is a polynomial and $\grad h$ is an automorphism. If $\grad h$ is Pascal finite of depth $m$, then $\Delta^k_{\grad h}(Q)=0$ for every $Q\in\kk[X]$ and every $k>(m-1)\deg Q$.
\end{corollary}
\begin{proof}
The first statement follows from Proposition~\ref{prop:PL} and Theorem~\ref{thm:tower}. For the second, put $F=\grad h$ and $\sigma(Q)=Q\circ F$, so that $\Delta_F=\sigma-1$. Since $\Delta_F^m(X)=0$, the iterates $F^{\circ j}=\sigma^j(X)=\sum_{l<m}\binom jl\Delta_F^l(X)$, $j\ge0$, are polynomials in $j$ of degree at most $m-1$ with coefficients in $\kk[X]$. Hence $\sigma^j(Q)=Q\circ F^{\circ j}$ is a polynomial in $j$ of degree at most $(m-1)\deg Q$, and $\Delta_F^k(Q)=\sum_{i=0}^k(-1)^{k-i}\binom ki\sigma^i(Q)$ is its $k$-th forward difference at $j=0$, which vanishes for $k>(m-1)\deg Q$.
\end{proof}

\begin{corollary}\label{cor:SN}
Let $h=\frac12\ip{X,X}+f$ be a Keller potential in normal form and suppose that $\Hess f$ is strongly nilpotent. Then $\grad h=X+\grad f$ is Pascal finite. If $m$ is its Pascal depth, then $\deg(\grad h)^{-1}\le(\deg\grad h)^{m-1}$ and $h^{*}=\sum_{k\le(m-1)\deg g}(-1)^k\Delta_{\grad h}^k(g)$.
\end{corollary}
\begin{proof}
$J(\grad f)=\Hess f$, so $\grad h=X+\grad f$ is a map with strongly nilpotent $J(\grad f)$. Such maps are Pascal finite by \cite[Theorem~5.1]{AH27}, and the degree bound holds for every Pascal finite map of the form $X+H$ by \cite[Proposition~5.3]{AH27}. Apply Proposition~\ref{prop:PL} and Corollary~\ref{cor:single-pascal}.
\end{proof}

\begin{lemma}\label{lem:qt}
Let $F=X+H$ over an arbitrary field. Then $F^{-1}=X-H$ if and only if $\Delta_F^2(X)=0$.
\end{lemma}
\begin{proof}
$(X-H)\circ F=X+H-H\circ F=X-\Delta_F(H)=X-\Delta_F^2(X)$. If this equals $X$, then $\sigma_F(P)=P\circ F$ satisfies $\sigma_F\circ\sigma_{X-H}=\mathrm{id}$, so $\sigma_F$ is a surjective endomorphism of the Noetherian ring $\kk[X]$, hence bijective, since a surjective endomorphism of a Noetherian ring is injective, and $X-H$ is a two-sided inverse of $F$.
\end{proof}

\begin{remark}\label{rem:isotropic}
Let $\alpha_1,\dots,\alpha_r\in\kk^n$ satisfy $\ip{\alpha_i,\alpha_j}=0$ for all $i,j$, let $q\in\kk[y_1,\dots,y_r]$ have order at least $3$, and let $f=q(\ip{\alpha_1,X},\dots,\ip{\alpha_r,X})$ and $h=\frac12\ip{X,X}+f$. Then
\[
\grad f=\sum_{i=1}^r(\partial_iq)\bigl(\ip{\alpha_1,X},\dots,\ip{\alpha_r,X}\bigr)\,\alpha_i ,
\]
hence $\ip{\alpha_j,X+\grad f}=\ip{\alpha_j,X}$ for all $j$, so that $\grad f\circ(X+\grad f)=\grad f$ and $(X-\grad f)\circ(X+\grad f)=X$. Moreover $\Hess f$ has square zero, so $h$ is a Keller potential in normal form whose gradient is a quasi-translation, of Pascal depth at most $2$ by Lemma~\ref{lem:qt}. By \cite[Theorem~1.2(i)]{vdEW08}, if $f=\sum_{j=1}^s\ip{\alpha_j,X}^d$ with $d\ge4$, pairwise linearly independent isotropic $\alpha_j$ and $\dim\langle\alpha_1,\dots,\alpha_s\rangle\in\{1,2,s-1,s\}$ is Hessian nilpotent, then $\ip{\alpha_i,\alpha_j}=0$ for all $i,j$, so that $f$ is of the form considered above with $q=\sum_jy_j^d$.
\end{remark}

\begin{theorem}[dimension at most three]\label{thm:dim3}
Let $n\le3$ and let $h$ be a Keller potential in normal form. Then $\grad h=X+H$ with $JH$ strongly nilpotent. In particular, $\grad h$ is Pascal finite and Corollary~\ref{cor:SN} applies.
\end{theorem}
\begin{proof}
If $\deg h\le2$, then $H=0$. Let $\deg h\ge3$. By \cite[Theorem~2.1]{dB15} there is $T\in\mathrm{GL}_n(\kk)$ such that all entries of $\Hess(h\circ T)$ below the anti-diagonal vanish. Put $F=\grad(h\circ T)=T^{\mathsf T}(\grad h)\circ T$ and $Q_i=F_{n+1-i}$. The entries of row $n+1-i$ of $JF$ in the columns $j>i$ vanish, so $Q_i\in\kk[x_1,\dots,x_i]$. Reversing the order of the rows of $\Hess(h\circ T)$ gives a lower triangular matrix with diagonal entries $\partial_iQ_i$. Hence these entries are polynomials whose product is $\pm\det\Hess(h\circ T)=\pm(\det T)^2\det\Hess h\in\kk^\times$, Since the only units of $\kk[X]$ are the nonzero constants, $\partial_iQ_i=c_i\in\kk^\times$. As $Q_i\in\kk[x_1,\dots,x_i]$ and $\partial_iQ_i=c_i$, the polynomial $Q_i-c_ix_i$ does not depend on $x_i$, so $Q_i-c_ix_i\in\kk[x_1,\dots,x_{i-1}]$, compare \cite[Lemma~2.5]{dB15}. Let $P$ be the anti-diagonal permutation matrix, so that $Q=PF$, $Q(0)=0$ and $M:=JQ(0)=PT^{\mathsf T}T$ is lower triangular with diagonal $(c_1,\dots,c_n)$. Since $M^{-1}$ is lower triangular with diagonal $(c_1^{-1},\dots,c_n^{-1})$, the map $M^{-1}Q=X+\tilde H$ satisfies $\tilde H_i\in\kk[x_1,\dots,x_{i-1}]$, so $J\tilde H$ is strictly lower triangular and strongly nilpotent. Finally
\[
M^{-1}Q=(PT^{\mathsf T}T)^{-1}PT^{\mathsf T}(\grad h)\circ T=T^{-1}\circ\grad h\circ T,
\]
so $\tilde H=T^{-1}H\circ T$ and $JH(Y)=T\,J\tilde H(T^{-1}Y)\,T^{-1}$. Products $JH(Y^{(1)})\cdots JH(Y^{(n)})$ are therefore conjugate to products of strictly lower triangular matrices and vanish. Pascal finiteness follows from \cite[Theorem~5.1]{AH27}.
\end{proof}

\begin{remark}\label{rem:dim2}
For $n=2$ the proof gives $\tilde H=(0,q(x_1))$, so that $\grad h$ is conjugate to a quasi-translation and $\Delta_{\grad h}^2(X)=0$. Theorem~\ref{thm:dim3} rests on \cite[Theorem~2.1]{dB15}, which fails for $n=4$ by the example \cite[(3)]{dB15}. Independently of this, Proposition~\ref{prop:phi0} below shows that the conclusion of Theorem~\ref{thm:dim3} fails for $n=4$. The invertibility of gradient Keller maps for $n\le3$ is the theorem of de Bondt \cite{dB15}. Theorem~\ref{thm:dim3} adds the structural statement that the normal form has strongly nilpotent nonlinear part, and with it Pascal finiteness and Corollary~\ref{cor:SN}.
\end{remark}

\begin{proposition}[degree of the inverse in dimension three]\label{prop:deg3}
Let $n=3$ and let $h$ be a Keller potential of degree $d\ge2$. Then $\grad h$ is an automorphism, $\deg(\grad h)^{-1}\le2d-3$, and $\delta^\alpha(g)=0$ for all $\alpha$ with $|\alpha|=2d-1$. For $h=x_1x_3+\frac12x_2^2+x_1^{d-1}x_2$ one has $\deg(\grad h)^{-1}=2d-3$, and $\delta^\alpha(g)\neq0$ for some $\alpha$ with $|\alpha|=2d-2$.
\end{proposition}
\begin{proof}
If $d=2$, then $\grad h$ is affine with invertible linear part and $2d-3=1$. Let $d\ge3$. As in the proof of Theorem~\ref{thm:dim3}, there is $T\in\mathrm{GL}_3(\kk)$ such that all entries of $\Hess(h\circ T)$ below the anti-diagonal vanish. Put $F=\grad(h\circ T)=T^{\mathsf T}(\grad h)\circ T$. Then $\deg(h\circ T)=d$, and $F$ is an automorphism if and only if $\grad h$ is, with inverses of the same degree. By the proof of Theorem~\ref{thm:dim3}, $F_3\in\kk[x_1]$, $F_2\in\kk[x_1,x_2]$, and $\partial_1F_3=c_1$, $\partial_2F_2=c_2$, $\partial_3F_1=c_3$ with $c_1,c_2,c_3\in\kk^\times$. The symmetry of $JF$ gives $c_3=\partial_3F_1=\partial_1F_3=c_1$ and $\partial_2F_1=\partial_1F_2$. Hence
\[
F_3=c_1x_1+a,\qquad F_2=c_2x_2+q(x_1),\qquad F_1=c_1x_3+q'(x_1)\,x_2+s(x_1),
\]
with $a\in\kk$ and $q,s\in\kk[x_1]$. Since $\deg F\le d-1$, one has $r:=\deg q\le d-1$, with $\deg0=-\infty$, and $\deg s\le d-1$. The equation $F(x)=y$ is solved by
\[
x_1=\frac{y_3-a}{c_1},\qquad x_2=\frac{y_2-q(x_1)}{c_2},\qquad x_3=\frac{y_1-q'(x_1)\,x_2-s(x_1)}{c_1},
\]
so $F$ is an automorphism, and the degrees of the coordinates of $F^{-1}$ in $y$ are at most $1$, $d-1$ and $\max(1,2r-1,d-1)\le2d-3$. By the proof of Theorem~\ref{thm:tower}, $\deg h^{*}\le\deg(\grad h)^{-1}+1\le2d-2$, so $\partial^\alpha h^{*}=0$ for $|\alpha|=2d-1$, and Proposition~\ref{prop:tower} gives $\delta^\alpha(g)=0$.

For $h=x_1x_3+\frac12x_2^2+x_1^{d-1}x_2$ one has $\det\Hess h=-1$ and
\[
\grad h=\bigl(x_3+(d-1)x_1^{d-2}x_2,\ x_2+x_1^{d-1},\ x_1\bigr),
\]
with inverse
\[
(\grad h)^{-1}(y)=\bigl(y_3,\ y_2-y_3^{d-1},\ y_1-(d-1)y_3^{d-2}y_2+(d-1)y_3^{2d-3}\bigr)
\]
of degree $2d-3$. Since $\grad h^{*}=(\grad h)^{-1}$ by Lemma~\ref{lem:leg}, the polynomial $h^{*}$ has degree $2d-2$. Hence $\partial^\alpha h^{*}$ is a nonzero constant for some $\alpha$ with $|\alpha|=2d-2$, and $\delta^\alpha(g)=(\partial^\alpha h^{*})\circ\grad h\neq0$ by Proposition~\ref{prop:tower}.
\end{proof}

\begin{remark}\label{rem:deg3}
For $n=3$ and $d\ge3$ one has $2d-3<(d-1)^2$. Thus the bound of \cite[Chapter~I, Corollary~1.4]{BCW82} is not attained by gradient automorphisms in three variables, and the order $N=(d-1)^2+2$ of Theorem~\ref{thm:tower} may be replaced by $2d-1$, which is optimal by the example of Proposition~\ref{prop:deg3}. The proof of Corollary~\ref{cor:pointwise} uses only an upper bound $B$ for $\deg G$, so for $n=3$ one may take $B=2d-3$ there. For $n=2$, Remark~\ref{rem:dim2} gives $\deg(\grad h)^{-1}\le d-1$. For $n=2$ and $n=3$ one therefore has $\deg(\grad h)^{-1}\le(n-1)(d-2)+1$, a bound linear in $d$. The proof for $n=3$ rests on \cite[Theorem~2.1]{dB15}, which fails for $n=4$. For the automorphism $\grad\phi_0$ of Proposition~\ref{prop:phi0} one has $d=4$ and $\deg(\grad\phi_0)^{-1}=3$.
\end{remark}

\begin{proposition}\label{prop:phi0}
Let $i\in\kk$ satisfy $i^2=-1$ and
\[
\phi_0=z_1z_4+z_2z_3+z_1(z_2^2+z_3^2)+2z_1^2z_2z_3+i\,z_1^2(z_2^2-z_3^2)\in\kk[z_1,\dots,z_4].
\]
Then $\det\Hess\phi_0=1$, the map $\grad\phi_0$ is a polynomial automorphism with inverse $\grad\psi_0$,
\[
\psi_0=w_1w_4+w_2w_3-w_4(w_2^2+w_3^2)+w_4^2\bigl(i\,w_2^2+2w_2w_3-i\,w_3^2\bigr),
\]
and the normal form of $\grad\phi_0$ at $0$ is not Pascal finite.
\end{proposition}
\begin{proof}
Put $u=(z_2,z_3)^{\mathsf T}$ and
\[
S(z_1)=\begin{pmatrix}2z_1+2iz_1^2&1+2z_1^2\\ 1+2z_1^2&2z_1-2iz_1^2\end{pmatrix},\qquad \det S(z_1)=-1 .
\]
Then $\phi_0=z_1z_4+\frac12u^{\mathsf T}S(z_1)u$, so that
\[
\grad\phi_0=\bigl(z_4+\tfrac12u^{\mathsf T}S'(z_1)u,\ S(z_1)u,\ z_1\bigr).
\]
For $p=(w_1,w_2,w_3,w_4)\in\kk^4$ put $w=(w_2,w_3)^{\mathsf T}$. The equation $\grad\phi_0(z)=p$ is equivalent to
\[
z_1=w_4,\qquad u=S(w_4)^{-1}w=-\adj S(w_4)\,w,\qquad z_4=w_1-\tfrac12u^{\mathsf T}S'(w_4)u .
\]
Hence $\grad\phi_0$ is bijective, and these formulas define its polynomial inverse $G$. Moreover
\begin{align*}
\ip{p,G(p)}-\phi_0(G(p))&=w_1w_4+w^{\mathsf T}u+w_4z_4-w_4z_4-\tfrac12u^{\mathsf T}S(w_4)u\\
&=w_1w_4+\tfrac12w^{\mathsf T}S(w_4)^{-1}w .
\end{align*}
Since $S^{-1}=-\adj S$, the last term equals $-\frac12\bigl((2w_4-2iw_4^2)w_2^2-2(1+2w_4^2)w_2w_3+(2w_4+2iw_4^2)w_3^2\bigr)$. Hence the Legendre transform of $\phi_0$ is $\psi_0$, and $G=\grad\psi_0$ by Lemma~\ref{lem:leg}. The variable $z_4$ occurs only in $z_1z_4$, so the Laplace expansion of $\Hess\phi_0$ along the row and the column of $z_4$ gives $\det\Hess\phi_0=-\det S(z_1)=1$.

Since $\phi_0(0)=0$, $\grad\phi_0(0)=0$ and $\Hess\phi_0(0)=L$, the antidiagonal permutation matrix, with $L^{-1}=L$, a normal form of $\grad\phi_0$ at $0$ is $A^{-1}\circ F\circ A$ with $F=L\grad\phi_0$ and any $A\in\mathrm{GL}_4(\kk)$ with $A^{\mathsf T}LA=I$, by Remark~\ref{rem:normalform}. Linear conjugation preserves Pascal finiteness and the degrees of iterates, so it suffices to consider $F$. If the normal form were Pascal finite of depth $m$, then its iterates $\sum_{l<m}\binom jl\Delta^l(X)$, and hence the iterates $F^{\circ j}$, would have bounded degree. Now $F(z)=\bigl(z_1,M(z_1)u,z_4+\frac12u^{\mathsf T}S'(z_1)u\bigr)$ with
\[
M(z_1)=\begin{pmatrix}0&1\\1&0\end{pmatrix}S(z_1)=I+z_1N_1+z_1^2N_2,
\]
where $N_1=\begin{pmatrix}0&2\\2&0\end{pmatrix}$ and $N_2=\begin{pmatrix}2&-2i\\2i&2\end{pmatrix}$. One has $N_2^2=4N_2\neq0$. Since $F$ fixes $z_1$, the $u$-part of $F^{\circ k}$ is $M(z_1)^ku$. The matrix $M(z_1)^k$ has degree at most $2k$ in $z_1$, and its coefficient of $z_1^{2k}$ is $N_2^k=4^{k-1}N_2\neq0$, since only the products of the terms $z_1^2N_2$ contribute to this coefficient. Hence $M(z_1)^ku$ contains the nonzero part $z_1^{2k}N_2^ku$ of degree $2k+1$, so $\deg F^{\circ k}\ge2k+1$ is unbounded, and the normal form is not Pascal finite.
\end{proof}

\begin{remark}\label{rem:phi0}
Proposition~\ref{prop:phi0} shows that Pascal finiteness of the gradient is a sufficient and not a necessary condition in $\HC_4$, and that the statement ``every gradient Keller map in four variables in normal form is Pascal finite'' is false. For $F=X+H$ as in the proof one has $\operatorname{tr}JH=4z_1^2$, so $JH$ is not nilpotent, since a nilpotent matrix has trace zero, hence not strongly nilpotent, and \cite[Theorem~5.1]{AH27} does not apply. The map $F$ is triangular. It fixes $z_1$, its middle coordinates $M(z_1)u$ are linear in $u$ over $\kk[z_1]$ with $\det M=1$, so that $u$ is recovered from them, and $z_4$ is then recovered from the last coordinate. This is the mechanism behind the invertibility.
\end{remark}

\begin{proposition}[the Pascal--Legendre series of $\phi_0$]\label{prop:phi0-g}
Let $\tilde h=\phi_0\circ A$ be the normal form of the proof of Proposition~\ref{prop:phi0} and $\tilde g$ its Legendre potential. Then $\Delta^m_{\grad\tilde h}(\tilde g)\neq0$ for every $m$. In particular, the sufficient condition of Corollary~\ref{cor:single-pascal} fails for $\phi_0$.
\end{proposition}
\begin{proof}
We use the notation of the proof of Proposition~\ref{prop:phi0}. One has $\tilde g=g\circ A$ with $g=\ip{z,\grad\phi_0}-\phi_0$, and $\grad\tilde h=A^{-1}\circ F\circ A$, so $\tilde g\circ(\grad\tilde h)^{\circ j}=(g\circ F^{\circ j})\circ A$. If $\Delta^m(\tilde g)=0$, then, as in the proof of Corollary~\ref{cor:single-pascal}, the polynomials $\tilde g\circ(\grad\tilde h)^{\circ j}=\sum_{l<m}\binom jl\Delta^l(\tilde g)$, $j\ge0$, have bounded degree, and so have the polynomials $g\circ F^{\circ j}$. We show that $\deg(g\circ F^{\circ k})\ge4k+4$.

By induction on $k$,
\[
F^{\circ k}=\bigl(z_1,\ M(z_1)^ku,\ z_4+c_k\bigr),\qquad c_k=\tfrac12\sum_{j=0}^{k-1}\bigl(M(z_1)^ju\bigr)^{\mathsf T}S'(z_1)\,M(z_1)^ju .
\]
Since $M(z_1)^j$ has degree at most $2j$ in $z_1$ and $S'(z_1)$ has degree one, $c_k$ has degree at most $4k-3$ in $z_1$. A direct computation gives
\[
g=z_1z_4+z_2z_3+2z_1(z_2^2+z_3^2)+3z_1^2\,u^{\mathsf T}\Theta u,\qquad \Theta=\begin{pmatrix}i&1\\1&-i\end{pmatrix}.
\]
Put $v=(v_1,v_2)^{\mathsf T}=M(z_1)^ku$. Then
\[
g\circ F^{\circ k}=z_1(z_4+c_k)+v_1v_2+2z_1(v_1^2+v_2^2)+3z_1^2\,v^{\mathsf T}\Theta v .
\]
Since $v$ has degree at most $2k$ in $z_1$, the first term has degree at most $4k-2$ in $z_1$ and the next two terms at most $4k+1$. The coefficient of $z_1^{2k}$ in $v$ is $N_2^ku$, so the coefficient of $z_1^{4k+2}$ in $g\circ F^{\circ k}$ is $3\,u^{\mathsf T}(N_2^k)^{\mathsf T}\Theta N_2^ku$. Now $\Theta N_2=4\Theta$, and since $\Theta$ is symmetric also $N_2^{\mathsf T}\Theta=4\Theta$. Hence $(N_2^k)^{\mathsf T}\Theta N_2^k=16^k\Theta$, and the coefficient equals $3\cdot16^k(i\,z_2^2+2z_2z_3-i\,z_3^2)\neq0$. Therefore $\deg(g\circ F^{\circ k})\ge4k+4$.
\end{proof}

\begin{remark}\label{rem:phi0-SS}
By Theorem~\ref{thm:member} and Remark~\ref{rem:SS}, for $h=\phi_0$ the remainder of $g$ on division by a Gr\"obner basis of $(\partial_1\phi_0-T_1,\dots,\partial_4\phi_0-T_4)$ for a block order with $z\gg T$ is $r=\psi_0(T)$. Thus one Gr\"obner basis and one reduction determine the inverse $\grad\psi_0$ of $\grad\phi_0$, whereas the normal form of $\grad\phi_0$ is not Pascal finite.
\end{remark}

\begin{remark}\label{rem:symmetry}
Write $G=\sum_sG_s$ for the homogeneous parts of the formal inverse. Comparing parts of degree $s\ge2$ in $G\circ F=X$ gives $G_s=-\sum_{j<s}[G_j\circ F]_s$. Hence, if $B\ge1$ and $G_s=0$ for $B<s<s_0$, then $G_{s_0}=-[G_{\le B}\circ F]_{s_0}$, and the vanishing of $G_s$ for $B<s\le eB$, $e=\deg F$, is equivalent to the single identity $G_{\le B}\circ F=X$. An induction on the degree is therefore equivalent to the symmetry condition of \cite[Corollary~6, Remark~8]{Ad26} on the polynomial side, and it requires an argument from outside the algorithm.
\end{remark}

\section{Limits of the criteria}\label{sec:scope}

\begin{proposition}\label{prop:scope}
Let $u=1+x_1x_2$,
\begin{gather*}
A=y_1u^3+3x_1y_2u^2-x_1^3y_3,\\
B=y_1x_2^2u(4+3x_1x_2)+y_2\bigl(x_2+3x_1x_2^2(4+3x_1x_2)\bigr)+y_3(2x_1-3x_1^2x_2),
\end{gather*}
and $\Psi=A^2+13A+2B\in\mathbb{Z}[x_1,x_2,y_1,y_2,y_3]$. Then $\deg\Psi=14$, $\det\Hess\Psi=128$, and $\grad\Psi(P_+)=\grad\Psi(P_-)=(0,0,-\frac12,0,0)$ for $P_\pm=(\pm1,\mp\frac32,0,0,0)$ \cite[Theorem~6.1]{MY26}. Consequently, for $h=\Psi$ and every $N\ge1$ there is a multi-index $\alpha$ with $|\alpha|=N$ and $\delta^\alpha(g)\neq0$, and for a normal form $\tilde\Psi$ of $\Psi$ and its Legendre potential $\tilde g$ one has $\Delta_{\grad\tilde\Psi}^k(\tilde g)\neq0$ for every $k$, by Corollary~\ref{cor:single-pascal}. In particular, the symmetry of the Hessian, the constancy of the Hessian determinant and the identity $\delta_i(g)=X_i$ do not by themselves imply that the formal Legendre transform is a polynomial.
\end{proposition}
\begin{proof}
Since $\grad\Psi$ is not injective, it has no inverse, and by Theorem~\ref{thm:tower}, $\Psi^{*}$ is not a polynomial. Apply Theorem~\ref{thm:tower} in the form valid for every order, and Corollary~\ref{cor:single-pascal}.
\end{proof}

\begin{remark}\label{rem:scope}
The potential $\Psi$ has all properties of a Keller potential used in Sections~\ref{sec:nambu}, \ref{sec:criterion}, \ref{sec:gradient} and~\ref{sec:pascal}. The matrix $J\grad\Psi=\Hess\Psi$ is symmetric with constant determinant, the Nambu derivations of $\grad\Psi$ are derivations of $\mathbb{Q}[x,y]$, and $\delta_i(g_\Psi)$ is the $i$-th coordinate function for the Legendre potential $g_\Psi$. The criteria of the paper recognize the non-invertibility of $\grad\Psi$. By Theorems~\ref{thm:tower} and~\ref{thm:member}, for every $N\ge1$ some $\delta^\alpha(g_\Psi)$ with $|\alpha|=N$ is nonzero, and $g_\Psi\notin\kk[\partial_1\Psi,\dots,\partial_5\Psi]$. Since a Keller map which is birational is an automorphism \cite[Chapter~I, Theorem~2.1]{BCW82}, one has $\mu(\grad\Psi)\ge2$, and by Theorem~\ref{thm:primitive} the minimal polynomial of $g_\Psi$ over $\kk(\grad\Psi)$ has degree $\mu(\grad\Psi)\ge2$. Moreover, since $\Psi^{*}$ is not a polynomial, \cite[Theorem~1.3]{Men06} shows that $\Hess\Psi$ is indefinite at every point of $\mathbb R^5$. Over $\mathbb R$ every counterexample to $\HC_n$ has this property.
\end{remark}

\section*{Acknowledgements}
This research was supported by the AGH University of Krakow within a subsidy of the Polish Ministry of Science and Higher Education (grant no.~16.16.420.054).

\end{document}